# KUMMER THEORY OF DIVISION POINTS OVER DRINFELD MODULES OF RANK ONE

WEN-CHEN CHI[†] AND ANLY LI[‡]

ABSTRACT. A Kummer theory of division points over rank one Drinfeld $A = \mathbb{F}_q[T]$- modules defined over global function fields was given. The results are in complete analogy with the classical Kummer theory of division points over the multiplicative algebraic group $\mathbb{G}_m$ defined over number fields.

## 0. INTRODUCTION

Let $K$ be a number field and let $\bar{K}$ be a fixed algebraic closure of $K$. For any positive integer $n$, let $\mu_n$ be the group of $n$-th roots of unity in $\bar{K}$. Let $G(n) = Gal(K(\mu_n)/K)$. For $K = \mathbb{Q}$, $G(n) \cong (\mathbb{Z}/n\mathbb{Z})^*$, and for any number field $K$, $G(l) \cong (\mathbb{Z}/l\mathbb{Z})^*$ for almost all prime numbers $l$.

For a finitely generated multiplicative subgroup $\Gamma$ of rank $r$ in $\mathbb{G}_m(K) = K^*$, $\Gamma$ is of finite index in its division group $\Gamma'$ in $K^*$. One considers the tower of Kummer extensions $K \subset K(\mu_n) \subset K(\mu_n, \Gamma^{\frac{1}{n}})$, where $K(\mu_n, \Gamma^{\frac{1}{n}})$ is the Galois extension of $K$ by adjoining the $n$-th roots of unity and the $n$-division points of $\Gamma$ in $\mathbb{G}_m(\bar{K})$. Let $H_\Gamma(n) = Gal(K(\mu_n, \Gamma^{\frac{1}{n}})/K(\mu_n))$ and $G_\Gamma(l) = Gal(K(\mu_n, \Gamma^{\frac{1}{n}})/K)$.

Classical Kummer theory of division points over the multiplicative algebraic group $\mathbb{G}_m$ over $K$ asserts the following well-known results (see [10, Theorem 4.1] ):

(i). For $K = \mathbb{Q}$, if $n$ is prime to $2[\Gamma' : \Gamma]$, then $H_\Gamma(n)$ is isomorphic to the direct product of $r$ copies of the abelian group $\mu_n$.

(ii). For any number field $K$, $H_\Gamma(l) \cong \mu_l \times ... \times \mu_l$ ($r$-copies) for almost all prime numbers $l$.

In this paper, we provide an analogous Kummer theory over the additive algebraic group $\mathbb{G}_a$ with additional module structure in the function field setting. More precisely, let $k = \mathbb{F}_q(T)$ be the rational function field over a finite field $\mathbb{F}_q$ and let $L$ be a finite extension of $k$ in a fixed algebraic closure $\bar{k}$ of $k$. Let $A = \mathbb{F}_q[T]$ and let $\phi$ be a Drinfeld $A$-module of rank one defined over $L$, where $L$ is viewed as an $A$-field of generic characteristic ( see [5, section 4.4], for general definition of Drinfeld modules). In particular, the Carlitz module ( see [5, Chapter 3] ) is a rank one Drinfeld module defined over $k$.

For a monic polynomial $M$ in $A$, let $\Lambda_M^\phi$ be the $M$-torsion points of the Drinfeld module $\phi$. Explicitly, $\Lambda_M^\phi = \{\alpha \in \bar{k} \mid \phi_M(\alpha) = 0\}$, where $\phi_M(\alpha)$ denotes the action of $M$ on $\alpha$. It is known that for the Carlitz module, $Gal(k(\Lambda_M^\phi)/k) \cong (A/MA)^*$, and for any rank one Drinfeld module $\phi$ defined over an $A$-field $L$ of generic characteristic, we have $Gal(L(\Lambda_l^\phi)/L) \cong (A/lA)^*$ for almost all monic irreducible polynomials $l$ in $A$ (see [5, Theorem 7.7] ).

Let $\Gamma$ be a finitely generated $A$-submodule of rank $r$ in the additive group $(L,+)$. Let $\frac{1}{M}\Gamma = \{\alpha \in \bar{k} \mid \phi_M(\alpha) \in \Gamma\}$ be the $M$-division module of $\Gamma$ in $(\bar{k},+)$. Then we have the tower of Kummer extensions $L \subset L(\Lambda_M^\phi) \subset L(\Lambda_M^\phi, \frac{1}{M}\Gamma)$. Let $H_\Gamma(M) = Gal(L(\Lambda_M^\phi, \frac{1}{M}\Gamma)/L(\Lambda_M^\phi))$ and $G_\Gamma(M) = Gal(L(\Lambda_M^\phi, \frac{1}{M}\Gamma)/L)$. Analogous to the classical case, we have the following results:

(i). For the Carlitz module, except the case that $q = 2$ and $T|M$ or $T+1|M$, under some mild condition, we have $H_\Gamma(M) \cong \Lambda_M^\phi \times ... \times \Lambda_M^\phi$ ($r$-copies) as $A$-modules.

(ii). For any rank one Drinfeld module, we have that for almost all monic irreducible polynomials $l$, $H_\Gamma(l) \cong \Lambda_l^\phi \times ... \times \Lambda_l^\phi$ ($r$-copies) as $A$-modules.

Here, one of the main idea is to show that an $A$-module structure is naturally equipped on $H_\Gamma(M)$ ( or particularly, on $H_\Gamma(l)$). Then, for the Carlitz module case, a Kummer theory along the line of the classical theory ( see [9], Ch.VI, Section 11) can be developed with this $A$-module structure naturally equipped throughout the whole theory. This establishes the above result (i).

For a general rank one Drinfeld $A$-module, using the $A$-module structure of $H_\Gamma(l)$ together with the independence property given by L.Denis ( see [4], Theorem 5), the above result (ii) can be established easily. Our proof is essentially the same as that of Denis in [4] except the above $A$-module is naturally equipped throughout the whole theory.

As for more general afiine rings $\mathbb{A}$, the same result as (ii) should follow by the same proof provided a proof of the independence property of Denis were given accordingly for $\mathbb{A}$. This includes a revised canonical height function and Dirichlet lemma for general affine ring $\mathbb{A}$ ( see [4], Section 4).

## 1. Some basic properties of Drinfeld modules of rank one

Throughout this paper, let $k = \mathbb{F}_q(T)$ be the rational function field of one variable over a finite field $\mathbb{F}_q$ of $q$ elements, where $q = p^m$ for some prime number $p$. Let $A = \mathbb{F}_q[T]$ be the polynomial ring over $\mathbb{F}_q$ which is the subring of those rational functions regular outside the place $\infty$ associated to $\frac{1}{T}$. Let $\bar{k}$ be a fixed algebraic closure of $k$. In this section, we will briefly review some definitions and basic properties of Drinfeld modules of rank one. For a general reference, we should refer to Chapter 3, 4 of [5] and [7].

First, recall that Carlitz makes $A$ act as a ring of endomorphisms on the additive group of $\bar{k}$ as follows:

Let $\tau : \bar{k} \to \bar{k}$ be the Frobenius automorphism defined by $\tau(\alpha) = \alpha^q$ and let $\mu_T$ be the map defined by $\mu_T(\alpha) = T\alpha$. The substitution $T \mapsto \tau + \mu_T$ yields a ring homomorphism from $A$ into the $\mathbb{F}_q$-algebra $End(\bar{k})$ of all $\mathbb{F}_q$-endomorphisms of the additive group of $\bar{k}$. This provides $\bar{k}$ with the structure of an $A$-module which is called the Carlitz module .

Write $\alpha^M$ for the action of $M \in A$ on $\alpha \in \bar{k}$, then we have $\alpha^M = M(\tau + \mu_T)(\alpha)$. In particular, for $a \in \mathbb{F}_q$, $\alpha^a = a\alpha$ for all $\alpha \in \bar{k}$. If $d = degM$, then $\alpha^M = \sum_{i=0}^d [{}^M_i] \cdot \alpha^{q^i}$, where each $[{}^M_i]$ is a polynomial in $A$ of degree $(d-i)q^i$ such that $[{}^M_0] = M$ and $[{}^M_d]$ is the leading coefficient of $M$. In [3, Equation 1.6], Carlitz gives an explicit formula for these polynomials.

For $M \neq 0$ in $A$, let $\Lambda_M = \{\alpha \in \bar{k} \mid \alpha^M = 0\}$. Then $\Lambda_M$ is an $A$-submodule of $\bar{k}$ which is called the module of $M$-torsion points of the Carlitz module. One has the following properties:

(1.1) $\Lambda_M$ is a vector space over $\mathbb{F}_q$ of dimension $d$, where $d = degM$.

(1.2) $\Lambda_M$ is a cyclic $A$-module with $\Phi(M)$ generators, where $\Phi(M)$ is the order of the group of units of $A/(M)$.

In fact, if $\lambda$ is a given generator and $B \in A$, then $\lambda^B$ is a generator if and only if $B$ and $M$ are relatively prime. Moreover, $(M)$ is equal to the annihilator of $\Lambda_M$. Hence $\Lambda_M$ is $A$-isomorphic to $A/(M)$.

(1.3) The $M$-torsion points $\Lambda_M$ generate a finite abelian extension, namely, the $M$-th cyclotomic function field $k(\Lambda_M)$ over $k$ such that $Gal(k(\Lambda_M)/k) \cong (A/(M))^*$. The isomorphism was given by $\sigma_B \mapsto B$, where $\sigma_B(\lambda) = \lambda^B$ for a given generator $\lambda$ of $\Lambda_M$ over $A$. In particular, $J = \{\sigma_a \mid a \in \mathbb{F}_q^*\}$ is a subgroup of $Gal(k(\Lambda_M)/k)$ which is known to be the inertia group of any infinite prime of $k(\Lambda_M)$(see also [6, Proposition 1.3] ).

**Remark:** Since the $A$-action is given by a polynomial over $k$, the action of $Gal(k(\Lambda_M)/k)$ on $\Lambda_M$ commutes with the $A$-action. So, $\sigma_B(\lambda) = \lambda^B$, for all $\lambda \in \Lambda_M$.

A field $L$ is said to be an $A$-field if there is a ring homomorphism $\iota : A \to L$. An $A$-field $L$ is said to be of generic characteristic if the kernel of $\iota$ is zero; otherwise, $L$ is said to be of finite characteristic $\wp$, where $\wp = ker(\iota)$. Let $L$ be a finite extension of $k$ which is viewed as an $A$-field of generic characteristic. Then $\phi$ is said to be a rank one Drinfeld $A$-module defined over $L$ if $\phi$ is a ring homomorphism from $A$ to $End(\bar{L})$ with $\phi_T(X) = TX + aX^q$, for $a \neq 0, a \in L$. For example, the Carlitz module is a rank one Drinfeld module over $k$.

Let $\phi : A \to L\{\tau\}$ be a Drinfeld $A$-module of rank one defined over a finite extension $L$ of $k$, where $L$ is viewed as an $A$-field of generic characteristic. Denote $\phi_m(\alpha)$ to be the action of $m \in A$ on $\alpha \in \bar{k}$. For $m \neq 0$ in $A$, let $\Lambda_m^\phi = \{\alpha \in \bar{k} \mid \phi_m(\alpha) = 0\}$. Then $\Lambda_m^\phi$ is an $A$-submodule of $\bar{k}$ which was called the module of $m$-torsion points of the Drinfeld module $\phi$. One has the following properties:

(1.4) $\Lambda_m^\phi$ is a vector space over $\mathbb{F}_q$ of dimension $d$, where $d = deg\, m$.

(1.5) $\Lambda_m^\phi$ is an $A$-module which is isomorphic to $\frac{A}{mA}$.

(1.6) For every monic irreducible polynomial $l$ in $A$ which satisfies the following conditions:

  (a). $\phi$ has good reduction at the primes of $L$ lying over $l$.
  
  (b). $l$ is unramified in $L_s/k$, where $L_s/k$ is the maximal separable subextension of $L/k$.
  
  One have that $L(\Lambda_l^\phi)/L$ is a finite abelian extension such that $Gal(L(\Lambda_l^\phi)/L) \cong (A/lA)^*$ ( see [5, Theorem 7.7.1]).

**Remarks:**

(1). Since the $A$-action is given by a polynomial over $L$, the action of $Gal(L(\Lambda_m^\phi)/L)$ on $\Lambda_m^\phi$ commutes with the $A$-action.

(2). Let $\phi$ be a Drinfeld module of rank one over an $A$-field $L$ of generic characteristic. For monic irreducible polynomial $l$ in $A$ which satisfies the above conditions in (1.6), $Gal(L(\Lambda_l^\phi)/L)$ consists of elements of the form: $\sigma = \sigma_a : \lambda \mapsto \phi_a(\lambda)$, where $\lambda$ is a generator of $\Lambda_l^\phi$ over $A$ and $a \in A$ runs over a set of representatives of $(A/lA)^*$.

## 2. THE KUMMER THEORY OVER THE CARLITZ MODULE

In this section, let $\Gamma$ be a finitely generated $A$-submodule of the additive group $(k, +)$. For a nonconstant polynomial $M$ in $A$, let $\frac{1}{M}\Gamma = \{\alpha \in \bar{k} | \alpha^M \in \Gamma\}$ be the $M$-division module of $\Gamma$. Denote by $K = k(\Lambda_M)$ and $k_{M,\Gamma} = k(\Lambda_M, \frac{1}{M}\Gamma)$. Analogous to the classical Kummer theory over $\mathbb{Q}$, we are interested in the following tower of Kummer extensions $k \subset K \subset k_{M,\Gamma}$ with associated Galois groups:

$$
\begin{array}{c}
k_{M,\Gamma} \\
H_\Gamma(M) \mid \\
K \ G_\Gamma(M) \\
G_\Gamma(M)/H_\Gamma(M) \simeq G(M) \mid \\
k
\end{array}
$$

By (1.3), the Galois group $G(M)$ is isomorphic to $(A/(M))^*$. The main goal is to show that under some mild conditions, $H_\Gamma(M)$ is as large as possible.

Given $z \in k$ and let $f_z(u) = u^M - z$, where $u^M = \sum_{i=0}^{d} [{}^M_i] u^{q^i}$ is the polynomial in $u$ which gives the Carlitz $A$-action on $\bar{k}$ as we have defined in Section 1. Then $f_z(u) \in k[u]$ and it is easy to see that $f_z(u)$ is a separable polynomial of degree $q^d$. The following properties are well-known (see [5], [7]):

(2.1)  $W = \{\alpha + \lambda \mid \lambda \in \Lambda_M\}$, where $\alpha$ is any fixed root of $f_z(u)$ in $\bar{k}$, form the complete set of all roots of $f_z(u)$ in $\bar{k}$.

(2.2)  The splitting field $k_{M,z}$ of $f_z(u)$ over $k$, is a finite abelian extension of $K$ such that $H_{M,z} = Gal(k_{M,z}/K)$ is naturally embedded into $\Lambda_M$ by $\psi \mapsto \lambda(\psi)$ if $\psi(\alpha) = \alpha + \lambda(\psi)$. More generally, for any given finitely generated $A$-submodule $\Gamma$ of $(k, +)$, the composite of all $k_{M,z}, z \in \Gamma$, $k_{M,\Gamma}$, is also an abelian extension of $K$.

For any given $z \in k$, by (2.2), the Galois group $H_{M,z}$ is isomorphic to a subgroup $H_M$ of $\Lambda_M$. Considering the tower of Galois extensions $k \subset K \subset k_{M,z}$, the Galois group $Gal(K/k)$ acts naturally on $H_{M,z} = Gal(k_{M,z}/K)$ by conjugation. Keeping the notations in (1.3) and (2.2), we may identify the Galois group $Gal(K/k)$ with $(A/(M))^*$. Then this action is explicitly given as follows:

**Proposition 2.1.** $\sigma_{\bar{B}} \cdot \psi_\lambda = \psi_{\lambda^B}$, for all $\bar{B} \in (A/(M))^*, \psi_\lambda \in H_{M,z}$; where $\sigma_{\bar{B}}$ and $\psi_\lambda$ are given by $\sigma_{\bar{B}}(\lambda) = \lambda^B$ and $\psi_\lambda(\alpha) = \alpha + \lambda$.

*Proof.* For any given $\sigma_{\bar{B}} \in Gal(K/k)$, let $\sigma \in Gal(k_{M,z}/k)$ be an extension of $\sigma_{\bar{B}}$. Then, for any given $\psi_\lambda \in Gal(k_{M,z}/K)$, we have $\sigma_{\bar{B}} \cdot \psi_\lambda = \sigma \circ \psi_\lambda \circ \sigma^{-1}$. Note that $\sigma^{-1}(\alpha) = \alpha + \lambda'$ for some $\lambda' \in \Lambda_M$. Consequently,

$$
\begin{aligned}
(\sigma_{\bar{B}} \cdot \psi_\lambda)(\alpha) &= (\sigma \circ \psi_\lambda \circ \sigma^{-1})(\alpha) \\
&= (\sigma \circ \psi_\lambda)(\alpha + \lambda') \\
&= \sigma(\alpha + \lambda' + \lambda) \\
&= \alpha + \lambda^B \\
&= \psi_{\lambda^B}(\alpha).
\end{aligned}
$$

This completes the proof. □

Now we extend the preceding natural action of $(A/(M))^*$ on $Gal(k_{M,z}/K)$ to an action of $A/(M)$ on $Gal(k_{M,z}/K)$ as follows:

Given $f \in A$, in the case $q \neq 2$ or $q = 2$ but $T(T+1) \nmid M$, we can write $f \,(mod\,M)$ as a finite sum $\sum f_i(mod\,M)$ such that $(f_i, M) = 1$ for all $i$. This can be done by Chinese Remainder Theorem as follows:

**Proposition 2.2.** *Let $M$ be a fixed nonzero element in $A = \mathbb{F}_q[T]$. In the case that $q \neq 2$ or $q = 2$ with $T(T+1) \nmid M$, for any $f \in A$, $f\,(mod\,M)$ can be written as a finite sum $\sum f_i(mod\,M)$ with $(f_i, M) = 1$ for all $i$. On the other hand, if $q = 2$ and $T(T+1)|M$, then such a decomposition for $f\,(mod M)$ does not always exist.*

*Proof.* First, we assume that $q \neq 2$ or $q = 2$ with $T(T+1) \nmid M$. Let $M = P_1^{n_1}...P_t^{n_t}$, where $P_i$; $i = 1, 2, ..., t$; are distinct irreducible polynomials in $A$. In particular, $(P_i, P_j) = 1$ whenever $i \neq j$. The assertion is trivial when $(f, M) = 1$. If $M|f$, then $f = (f-1) + 1$ gives a desired finite sum for $f\,(mod\,M)$. So, we may assume that $(f, M) \neq 1$ and $M \nmid f$.

Consider $f\,(mod\,P_i^{n_i})$ for each $i = 1, 2, ..., t$. Let $I \subseteq \{1, ..., t\}$ be the set of indices $i$ such that $f \equiv a_i\,(mod\,P_i^{n_i})$ with $a_i \not\equiv 0(mod\,P_i^{n_i})$ and let $J = \{j \,|\, 1 \leq j \leq t$ and $P_j^{n_j}$ divides $f\} = \{1, ..., t\} \backslash I$.

If $2 \nmid q$, then by Chinese Remainder Theorem, there exist $f_1$ and $f_2$ in $A$ such that
$$f_1 \equiv \begin{cases} a_i/2\,(mod\,P_i^{n_i}) & , \text{ for } i \in I \text{ with } P_i \nmid a_i, \\ a_i - 1\,(mod\,P_i^{n_i}) & , \text{ for } i \in I \text{ with } P_i|a_i, \\ 1\,(mod\,P_j^{n_j}) & , \text{ for } j \in J, \end{cases}$$
and
$$f_2 \equiv \begin{cases} a_i/2\,(mod\,P_i^{n_i}) & , \text{ for } i \in I \text{ with } P_i \nmid a_i, \\ 1\,(mod\,P_i^{n_i}) & , \text{ for } i \in I \text{ and } P_i|a_i, \\ -1\,(mod\,P_j^{n_j}) & , \text{ for } j \in J. \end{cases}$$
Then $f \equiv f_1 + f_2\,(mod\,M)$ with $(f_1, M) = (f_2, M) = 1$.

For $2|q$, we discuss the two possible cases as follows:

Case(i): $q = 2^s$, where $s \geq 2$.

For $i \in I$ with $P_i \nmid a_i$, $a_i\,(mod\,P_i)$ is a nonzero element of the finite field $A/(P_i)$ which has at least two distinct nonzero elements. So we can choose a polynomial $b_i \in A$ with $(b_i, P_i) = 1$ such that $a_i + b_i \not\equiv 0\,(mod\,P_i)$. For $i \in I$ with $P_i|a_i$, it's obvious that $(a_i + 1, P_i) = 1$. Thus, for each $i \in I$, there always exists $b_i \in A$ with $(b_i, P_i) = 1$ such that $(a_i + b_i, P_i) = 1$. Apply Chinese Remainder Theorem, there exist $f_1$ and $f_2$ such that
$$f_1 \equiv \begin{cases} a_i + b_i\,(mod\,P_i^{n_i}) & , \text{ for } i \in I, \\ 1\,(mod\,P_j^{n_j}) & , \text{ for } j \in J, \end{cases}$$
and
$$f_2 \equiv \begin{cases} -b_i\,(mod\,P_i^{n_i}) & , \text{ for } i \in I, \\ -1\,(mod\,P_j^{n_j}) & , \text{ for } j \in J. \end{cases}$$
Then $f \equiv f_1 + f_2\,(mod\,M)$, where $(f_1, M) = (f_2, M) = 1$.

Case(ii): $q = 2$ with $T(T+1) \nmid M$.

In this situation, $deg P_i \geq 2$ for all $i \in I$. In particular, the finite field $A/(P_i)$ has at least two distinct nonzero elements. The same argument as in Case(i) would give a desired decomposition for $f\,(mod\,M)$.

Finally, assume that $q = 2$ and $T(T+1)|M$. Take an $f \in A$ such that $T|f$ and $T+1 \nmid f$. Suppose $f \equiv f_1 + ... + f_n\,(mod\,M)$ with $(f_i, M) = 1$ for all $i$, $1 \leq i \leq n$.

Then $f(0) = f(1) = n$. But $T|f$ implies that $f(0) = 0$ and $T + 1 \nmid f$ implies that $f(1) = 1$, which is a contradiction. Similarly, for $f \in A$ with $T \nmid f$ and $T + 1|f$, $f \pmod{M}$ cannot have the decomposition. This completes the proof. □

Thus in the case that $q \neq 2$ or $q = 2$ but $T(T + 1) \nmid M$, we can define, for $\psi_\lambda \in H_{M,z}$ and $\bar{f} \in A/(M)$,

$$\bar{f} \cdot \psi_\lambda = \sum_i \sigma_{\bar{f}_i} \cdot \psi_\lambda = \sum_i \psi_{\lambda^{f_i}} = \psi_{\sum \lambda^{f_i}}.$$

It is easy to check this action is independent of the decomposition $f \equiv \sum f_i \pmod{M}$ by noting that $\sum \lambda^{f_i} = \lambda^{\sum f_i} = \lambda^f$ which is independent of the choice of the $f_i$.

Therefore, this action is well-defined. Composing with the canonical map from $A$ to $A/(M)$, we have an $A$-action on $H_{M,z}$.

By the same way, $Gal(K/k)$ acts naturally on $Gal(k_{M,\Gamma}/K)$ by conjugation. In particular, we have an $(A/(M))^*$-action on $Gal(k_{M,\Gamma}/K)$. Denote this action by $\sigma_f \cdot \tau$, for $\sigma_f \in Gal(K/k)$ and $\tau \in Gal(k_{M,\Gamma}/K)$. It is easy to check that $(\sigma_f \cdot \tau)|_{k_{M,z}} = \sigma_f \cdot (\tau|_{k_{M,z}})$ for each $z \in \Gamma$. Except for the case $q = 2$ and $T(T+1)|M$, for each $\bar{f} \in A/(M)$, write $f \equiv \sum f_i \pmod{M}$ with $(f_i, M) = 1$. Then we can define $\bar{f} \cdot \tau = \sum \sigma_{f_i} \cdot \tau$, and hence $(\bar{f} \cdot \tau)|_{k_{M,z}} = \bar{f} \cdot (\tau|_{k_{M,z}})$ for each $z \in \Gamma$. This gives an $A$-action on $Gal(k_{M,\Gamma}/K)$.

Notice that under the natural embedding $\psi_\lambda \mapsto \lambda$ by (2.2), $H_{M,z}$ is isomorphic to a subgroup $H_M$ of $\Lambda_M$. The above definition obviously gives that $f \cdot \psi_\lambda = \psi_{\lambda^f}$. In particular, if $\lambda \in H_M$, then so is $\lambda^f$ for all $f \in A$. Thus, $H_M$ is an $A$-submodule of $\Lambda_M$. To summarize the above discussion, we have the following:

**Proposition 2.3.** *Except for the case $q = 2$ and $T(T + 1)|M$, we have:*

(1). *The $A$-action defined as above gives an $A$-module structure on $H_{M,z}$ and consequently gives an $A$-module structure on $H_{M,\Gamma}$.*
(2). *$H_M$ is an $A$-submodule of $\Lambda_M$ and $H_{M,z}, H_M$ are isomorphic as $A$-modules. Consequently, $H_{M,z}$ and $H_{M,\Gamma}$ are $A$-modules of exponent $M$.*

*Proof.* First, it is easy to check that $H_{M,z}$ is an $A$-module under the above well-defined $A$-action as follows:

(i). $H_{M,z}$ is known to be an abelian group.
(ii). For $f \in A$ and for $\psi_{\lambda_1}, \psi_{\lambda_2}$ in $H_{M,z}$ with $\lambda_1, \lambda_2$ in $H_M \subseteq \Lambda_M$, by Proposition 2.1, $f \cdot (\psi_{\lambda_1} + \psi_{\lambda_2}) = f \cdot \psi_{\lambda_1 + \lambda_2} = \psi_{(\lambda_1 + \lambda_2)^f} = \psi_{\lambda_1^f} + \psi_{\lambda_2^f} = f \cdot \psi_{\lambda_1} + f \cdot \psi_{\lambda_2}$.
(iii). Let $f, g \in A$ and let $\lambda \in H_M$. Write $f \equiv \sum f_i \pmod{M}$, $g \equiv \sum g_j \pmod{M}$, with $(f_i, M) = (g_j, M) = 1$ for all $i, j$. Then

$$(fg) \cdot \psi_\lambda = (\sum_{i,j} \sigma_{\bar{f}_i \bar{g}_j}) \cdot \psi_\lambda = \sum_{i,j} \psi_{\lambda^{f_i g_j}} = \psi_{\lambda^{fg}}.$$

On the other hand,

$$f \cdot (g \cdot \psi_\lambda) = f \cdot \psi_{\lambda^g} = \psi_{\lambda^{gf}} = (fg) \cdot \psi_\lambda.$$

Moreover, $(f + g) \cdot \psi_\lambda = \psi_{\lambda^{f+g}} = \psi_{\lambda^f} + \psi_{\lambda^g} = f \cdot \psi_\lambda + g \cdot \psi_\lambda$.

Finally, by Proposition 2.1, it is clear that $H_{M,z}$ and $H_M$ are isomorphic as $A$-modules. Since $H_M$ is of exponent $M$, so are the Galois groups $H_{M,z}$ and $H_{M,\Gamma}$. This completes the proof. □

**Remark:** $Gal(k(\Lambda_M)/k)$ acts on $H_{M,z}$ by conjugation and acts on $\Lambda_M$ naturally. By Proposition 2.1, $\sigma_{\bar{B}} \cdot \psi_\lambda = \psi_{\lambda^B} = \psi_{\sigma_{\bar{B}}(\lambda)}$, so $H_{M,z}$ and $H_M$ are isomorphic as $Gal(k(\Lambda_M)/k)$-modules as well.

Recall that $\Lambda_M$ is a cyclic $A$-module. Consequently, $H_M$ is a cyclic $A$-module and hence $H_{M,z}$ is a cyclic $A$-module. This leads to the following general definitions. To fix notations, let $E, F$ be extensions of $k$ in $\bar{k}$.

**Definition 2.4.**
(1). *An abelian Galois extension $E/F$ is said to be $A$-abelian if its Galois group has an $A$-module structure. Denote it by $(E/F, \cdot_A)$ to specify the $A$-module structure.*
(2). *An $A$-abelian extension $(E/F, \cdot_A)$ is said to be $A$-cyclic if its Galois group is a cyclic $A$-module. In this case, if $Gal(E/F) \cong A/(M)$, where $M$ is a monic polynomial, then we say that the $A$-cyclic extension $E/F$ is of order $M$.*

**Definition 2.5.** *An $A$-abelian extension $(E/F, \cdot_A)$ is said to be of exponent $M$ if its Galois group $G$ is a $M$-torsion $A$-module, i.e., $M \cdot_A \sigma = 1$ for all $\sigma \in G$.*

**Example:** Let $K = k(\Lambda_M)$ and $z \in k - k^M$. With the $A$-module structure defined in Proposition 2.3, $k_{M,z}/K$ is an $A$-cyclic extension of order $N$ dividing $M$ and $k_{M,\Gamma}/K$ is an $A$-abelian extension of exponent $M$.

**Remark:** For any field extension $E/k$ in $\bar{k}$ and for any automorphism $\sigma$ of $E$ over $k$, by the formula given by Carlitz, we have

$$\sigma(\alpha^M) = \sigma(\sum_{i=0}^{d} [{}_i^M]\alpha^{q^i})$$
$$= \sum_{i=0}^{d} [{}_i^M]\sigma(\alpha)^{q^i}$$
$$= \sigma(\alpha)^M \text{ for all } M \in A.$$

In other words, $\sigma$ is an $A$-module automorphism of the $A$-module $(E, +)$.

**Proposition 2.6.** *Assume $\Lambda_M \subseteq F$. If $(E/F, \cdot_A)$ is an $A$-cyclic extension of order $M$, then there exists $\alpha \in E$ such that $E = F(\alpha)$ and $\alpha$ satisfies an equation $X^M - a = 0$ for some $a \in F$.*

*Proof.* By definition, $G = Gal(E/F)$ is isomorphic to $A/(M)$ as $A$-modules. On the other hand, $\Lambda_M$ is isomorphic to $A/(M)$ as $A$-modules. Thus we have an $A$-isomorphism $f : G \to \Lambda_M$. If $\sigma$ is a generator of $G$ over $A$, then $\lambda = f(\sigma)$ is a generator of $\Lambda_M$ over $A$. Moreover, $f(B \cdot_A \sigma) = \lambda^B$ for all $B \in A$.

Consider the map $f$. Since $G$ acts trivially on $\Lambda_M$, we may view $f$ as a 1-cocycle of $G$ with values in the additive group $(E, +)$. It is well-known that $H^1(G, E) = 0$ by normal basis theorem. Consequently, there exists $\alpha \in E$ such that $f(\tau) = \tau\alpha - \alpha$ for all $\tau \in G$. In particular, $\sigma\alpha = \alpha + \lambda$, where $\sigma$ is a fixed generator of $G$ over $A$ and $\lambda = f(\sigma)$. For any $\tau \in G, \tau = B \cdot_A \sigma$ for some $B \in A$. Hence we have $B \cdot_A \sigma(\alpha) = \alpha + \lambda^B$ for all $B \in A$. We conclude that $\{\alpha + \lambda | \lambda \in \Lambda_M\}$ are distinct conjugates of $\alpha$ over $F$. This implies that $[F(\alpha) : F] \geq |A/(M)|$. Since $[E : F] = |A/(M)|$, we must have $E = F(\alpha)$. Furthermore, $\sigma(\alpha^M) = (\sigma(\alpha))^M = (\alpha + \lambda)^M = \alpha^M$ and for all $B \in A$, $B \cdot_A \sigma(\alpha^M) = (B \cdot_A \sigma(\alpha))^M = (\alpha + \lambda^B)^M = \alpha^M$. Thus $\alpha^M \in F$ and we let $a = \alpha^M$. This proves the assertion. □

Recall that, by (1.3), $J = \{\sigma_a | a \in \mathbb{F}_q^*\}$ is a subgroup of $G(M)$, where $\sigma_a(\lambda) = a\lambda$. This gives the following result by a well-known theorem of Sah (see [10, Theorem 5.1]).

**Proposition 2.7.** *Except for the case that $q = 2$ and $T|M$ or $T+1|M$, we have $H^1(G(M), \Lambda_M) = 0$.*

*Proof.* For $q \neq 2$, there exist elements $a$ and $a - 1 \in \mathbb{F}_q^*$ such that $\lambda \mapsto \sigma_a\lambda - \lambda$ is an automorphism of $\Lambda_M$. For $q = 2$ but neither $T$ nor $T + 1$ divides $M$, $(A/(M))^*$ contains elements $f$ and $f + 1$ such that $\lambda \mapsto \lambda^f - \lambda$ is an automorphism of $\Lambda_M$. Then by Sah's theorem, we have $H^1(G(M), \Lambda_M) = 0$. □

For the rest of this section, we assume that $q \neq 2$ or $q = 2$ but neither $T$ nor $T + 1$ divides $M$.

For the finitely generated $A$-submodule $\Gamma$ of $(k, +)$, let $\Gamma' = \frac{1}{M}\Gamma \cap k$ and define the exponent $e(\Gamma'/\Gamma)$ to be the unique monic polynomial with smallest degree such that $\Gamma'^{e(\Gamma'/\Gamma)} \subseteq \Gamma$. It is easy to check that $e(\Gamma'/\Gamma)$ is well-defined.

For each $a \in \Gamma$, let $\alpha \in \bar{k}$ be a root of the polynomial $f_a(X) = X^M - a$. Let $\sigma \in H_{M,\Gamma}$. Then $\sigma\alpha = \alpha + \lambda_\sigma$ for some $\lambda_\sigma \in \Lambda_M$. The map $\sigma \mapsto \lambda_\sigma$ is obviously a homomorphism of $H_{M,\Gamma}$ into $\Lambda_M$. Write $\lambda_\sigma = \sigma\alpha - \alpha$. It is easy to see that $\lambda_\sigma$ is independent of the choice of the root $\alpha$ of $X^M - a$. We denote $\lambda_\sigma$ by $<\sigma, a>$. The map $(\sigma, a) \mapsto <\sigma, a>$ gives us a map $H_{M,\Gamma} \times \Gamma \to \Lambda_M$.

**Proposition 2.8.** *The map $H_{M,\Gamma} \times \Gamma \to \Lambda_M$ given by $(\sigma, a) \mapsto <\sigma, a>$ is $A$-bilinear, so that the kernel on the left is $\{1\}$ and the kernel on the right is $\Gamma \cap K^M$.*

*Proof.* If $a, b \in \Gamma$ and $\alpha^M = a, \beta^M = b$, then $(\alpha + \beta)^M = a + b$ and hence $<\sigma, a + b> = \sigma(\alpha + \beta) - (\alpha + \beta) = (\sigma(\alpha) - \alpha) + (\sigma(\beta) - \beta) = <\sigma, a> + <\sigma, b>$ for all $\sigma \in H_{M,\Gamma}$. On the other hand, let $\sigma, \tau \in H_{M,\Gamma}$ and $a \in \Gamma$. If $\alpha^M = a$, then $\sigma\tau(\alpha) = \sigma(\alpha + \lambda_\tau) = \alpha + \lambda_\sigma + \lambda_\tau$. Hence

$$<\sigma\tau, a> = <\sigma, a> + <\tau, a>.$$

Moreover, for each $B \in A$ and for each $\sigma \in H_{M,\Gamma}$, by the definition of $A$-action **.** on $H_{M,\Gamma}$ (see the discussion above Proposition 2.3), we have $B \cdot \sigma(\alpha) = \alpha + \lambda_\sigma^B$. In other words,

$$<B \cdot \sigma, a> = B \cdot \sigma(\alpha) - \alpha = \lambda_\sigma^B = <\sigma, a>^B.$$

On the other hand, if $\alpha^M = a$, then $(\alpha^B)^M = a^B$. Hence $<\sigma, a^B> = \sigma(\alpha^B) - \alpha^B = (\sigma(\alpha))^B - \alpha^B = (\alpha + <\sigma, a>)^B - \alpha^B = <\sigma, a>^B$. This proves that the map $(\sigma, a) \mapsto <\sigma, a>$ is an $A$-module bilinear map from $H_{M,\Gamma} \times \Gamma$ to $\Lambda_M$.

Suppose $\sigma \in H_{M,\Gamma}$ such that $<\sigma, a> = 0$ for all $a \in \Gamma$. Then for every generator $\alpha$ of $k_{M,\Gamma}$ such that $\alpha^M = a$, we have $\sigma\alpha = \alpha$. Hence $\sigma = 1$ and the kernel on the left is $\{1\}$.

On the other hand, let $a \in \Gamma$ be such that $<\sigma, a> = 0$ for all $\sigma \in H_{M,\Gamma}$. Let $\alpha \in \bar{k}$ be such that $\alpha^M = a$. Consider the subfield $k_{M,a} = K(\alpha)$ of $k_{M,\Gamma}$. If $\alpha \notin K$, then there exists an automorphism of $K(\alpha)$ over $K$ which is not the identity. Extend this automorphism to $k_{M,\Gamma}$ and call this extension $\sigma$. Then clearly $<\sigma, a> \neq 0$. Thus the kernel on the right is $\Gamma \cap K^M$. □

Consequently, we have an $A$-module homomorphism $\varphi : \Gamma \to Hom_A(H_\Gamma(M), \Lambda_M)$. More precisely, for each $a \in \Gamma$, we have an $A$-module homomorphism

$$\varphi_a : H_\Gamma(M) \to \Lambda_M \text{ defined by } \varphi_a(\sigma) = \sigma\alpha - \alpha,$$

where $\alpha^M = a$.

**Theorem 2.9.** *Let $e_M(\Gamma) = g.c.d.(e(\Gamma^{'}/\Gamma), M)$ and let $\Gamma_\varphi$ be the kernel of $\varphi$. Then $\Gamma_\varphi^{e_M(\Gamma)} \subseteq \Gamma^M$.*

*Proof.* Let $a \in \Gamma_\varphi$ and $\alpha^M = a$. For each $\sigma \in G_\Gamma(M)$, define $\lambda_\sigma = \sigma\alpha - \alpha$. Then $\{\lambda_\sigma\}$ is a 1-cocycle of $G_\Gamma(M)$ in $\Lambda_M$. Since $a \in \Gamma_\varphi, \sigma\alpha = \alpha$ for all $\sigma \in H_\Gamma(M)$; this cocycle depends only on the class of $\sigma$ modulo the subgroup $H_\Gamma(M)$ of $G_\Gamma(M)$. We may view $\lambda_\sigma$ as a 1-cocycle of $G(M)$ in $\Lambda_M$. By Proposition 2.7, there exists a $\lambda_0 \in \Lambda_M$ such that $\lambda_\sigma = \sigma\lambda_0 - \lambda_0$ for all $\sigma \in G_\Gamma(M)$. Thus $\sigma(\alpha - \lambda_0) = \alpha - \lambda_0$ for all $\sigma \in G_\Gamma(M)$.

In other words, $\alpha - \lambda_0 \in k$. Since both $\alpha$ and $\lambda_0$ are in $\frac{1}{M}\Gamma$, we have $\alpha - \lambda_0 \in \Gamma^{'}$. This proves that $a = (\alpha - \lambda_0)^M \in (\Gamma^{'})^M$ for all $a \in \Gamma_\varphi$ and hence $\Gamma_\varphi \subseteq (\Gamma^{'})^M$. Since $e_M(\Gamma) = f \cdot e(\Gamma^{'}/\Gamma) + g \cdot M$ for some $f, g \in A$, we have that $\Gamma_\varphi^{e_M(\Gamma)} \subseteq (\Gamma^{'M})^{e_M(\Gamma)} \subseteq \Gamma^M$. This completes the proof. □

**Corollary 2.10.** *If $e_M(\Gamma) = 1$, then $\Gamma^M = \Gamma \cap K^M = \Gamma \cap k^M$. In this case, the pairing $H_\Gamma(M) \times \Gamma/\Gamma^M \to \Lambda_M$ is nondegenerate. Consequently, we have an $A$-module ( resp. $A/(M)$-module) isomorphism*

$$\varphi : \Gamma/\Gamma^M \to Hom_A(H_\Gamma(M), \Lambda_M)(\text{ resp. } \varphi : \Gamma/\Gamma^M \to Hom_{A/(M)}(H_\Gamma(M), \Lambda_M)).$$

*Proof.* By Theorem 2.9, the right kernel of the pairing $H_\Gamma(M) \times \Gamma \to \Lambda_M$ is contained in $\Gamma^M$. In other words, we have that $(\Gamma \cap K^M) \subseteq \Gamma^M$. On the other hand, $\Gamma^M \subseteq (\Gamma \cap k^M) \subseteq (\Gamma \cap K^M)$. We conclude that $\Gamma^M = \Gamma \cap K^M = \Gamma \cap k^M$. In particular, the pairing $H_\Gamma(M) \times \Gamma/\Gamma^M \to \Lambda_M$ is nondegenerate. By duality of $A$-( resp. $A/(M)$-) modules, we have the isomorphisms as stated. □

**Corollary 2.11.** *If $e_M(\Gamma) = 1$ and $\Gamma$ is free of rank $r$ with basis $\{a_1, ..., a_r\}$, let $\varphi_i = \varphi_{a_i}$, then the map $H_\Gamma(M) \to \Lambda_M \times ... \times \Lambda_M$ ( r-copies) given by $\sigma \mapsto (\varphi_1(\sigma), ..., \varphi_r(\sigma))$ is an $A$-module ( resp. $A/(M)$-module) isomorphism.*

*Proof.* It is easy to see that the map $\sigma \mapsto (\varphi_1(\sigma), ..., \varphi_r(\sigma))$ is injective. On the other hand, Corollary 2.10 implies that $H_\Gamma(M)$ has order $|A/(M)|^r$, which is also the order of $\Lambda_M \times ... \times \Lambda_M$ ( r-copies). Hence it is surjective. This proves the assertion. □

**Remarks:**

(i). Let $\Gamma$ be a finitely generated $A$-submodule of $(k, +)$ of rank $r$. By general theory of modules over principal ideal rings (see [1, Ch.VII, §4]) and Theorem 1 of [11], $\Gamma$ is isomorphic to a direct sum of the form $A \bigoplus ... \bigoplus A$ or $A \bigoplus ... \bigoplus A \bigoplus A/(N)$, where $N$ is a nonzero polynomial. If $e_M(\Gamma) = 1$, by Corollary 2.11, we have a noncanonical $A/(M)$-module isomorphism between $\Gamma/\Gamma^M$ and $H_\Gamma(M)$. If in addition $\Gamma \cong A \bigoplus ... \bigoplus A$ , or, $\Gamma \cong A \bigoplus ... \bigoplus A \bigoplus A/(N)$ and $M$ is relatively prime to $N$, then $H_\Gamma(M)$ is isomorphic to $\Lambda_M \times ... \times \Lambda_M$ ( r-copies).

(ii). If the orders of $\Lambda_M$ and $G(M)$ are relatively prime, for example,

$$M = \prod_{P|M} P \text{ is a product of distinct irreducible polynomials } P,$$

then $H^2(G(M), \Lambda_M) = 0$ by Cor.(10.2) in [2]. In this case, the orders of $H_\Gamma(M)$ and $G(M)$ are also relatively prime, so $H^2(G(M), H_\Gamma(M)) = 1$, where $G(M)$ acts on $H_\Gamma(M)$ by conjugation. In particular, the exact sequence $1 \to H_\Gamma(M) \to G_\Gamma(M) \to G(M) \to 1$ is split and $G_\Gamma(M)$ is a semidirect product of $H_\Gamma(M)$ by $G(M)$ ( see [2, CH. IV] ).

## 3. The Kummer theory over rank one Drinfeld $\mathbb{F}_q[T]$-modules

In this section, we will consider general rank one Drinfeld $A$-modules, and the following discussion will be similar with that given in the previous section. The main difference is that the Galois group of the cyclotomic extension can be completely determined for any nonzero polynomial in the Carlitz module case, while in general rank one case, it can only be determined under some condition ( see [5], Theorem7.7.1). For the convenience of the readers, we will also give the sketch of the proof.

Let $\phi$ be a Drinfeld $A$-module of rank one defined over a finite extension $L$ of $k$, where $L$ is viewed as an $A$-field of generic characteristic. For simplicity, we denote $L(\Lambda_m^\phi)$ by $L_m$ for all $m \neq 0$ in $A$. By definition, it is clear that the additive group of $L_m$ is an $A$-submodule of $\bar{k}$.

Given $z \in L$ and let $f_z(u) = \phi_m(u) - z$. Then $f_z(u) \in L[u]$ and it is easy to see that $f_z(u)$ is a separable polynomial of degree $q^d$, where $d$ is the degree of $m$. Similar to the discussions in Section 2, we can consider the splitting field $L_{m,z}$ of $f_z(u)$ over $L$, say $L_m(\alpha)$, where $\alpha$ is any fixed root of $f_z(u)$ in $\bar{k}$. And we have that $L_{m,z}$ is a finite abelian extension of $L_m$ such that $H_{m,z} = Gal(L_{m,z}/L_m)$ is naturally embedded into $\Lambda_m^\phi$ by $\psi \mapsto \lambda(\psi)$ if $\psi(\alpha) = \alpha + \lambda(\psi)$. More generally, for a given $A$-submodule $\Gamma$ of $L$, let $L_{m,\Gamma}$ be the composite of all $L_{m,z}; z \in \Gamma$. Then $L_{m,\Gamma}$ is also an abelian extension of $L_m$.

Throughout the rest of this section, $l$ will denote a monic irreducible polynomial in $A$ satisfying the following conditions:

(a). $\phi$ has good reduction at the primes of $L$ lying over $l$.
(b). $l$ is unramified in $L_s/k$, where $L_s/k$ is the maximal separable subextension of $L/k$.

For any given $z \in L$, via the above embedding, the Galois group $H_{l,z}$ is isomorphic to a subgroup $H_l$ of $\Lambda_l^\phi$. Considering the tower of Galois extensions $L \subset L_l \subset L_{l,z}$, the Galois group $Gal(L_l/L)$ acts naturally on $H_{l,z} = Gal(L_{l,z}/L_l)$ by conjugation. By identifying the Galois group $Gal(L_l/L)$ with $(A/lA)^*$, this action can be explicitly computed as Proposition 2.1, we have $\sigma_{\bar{a}} \cdot \psi_\lambda = \psi_{\phi_a(\lambda)}$, for all $\bar{a} \in (A/lA)^*, \psi_\lambda \in H_{l,z}$; where $\sigma_{\bar{a}}$ and $\psi_\lambda$ are given by $\sigma_{\bar{a}}(\lambda) = \phi_a(\lambda)$ and $\psi_\lambda(\alpha) = \alpha + \lambda$.

As in Section 2, we can extend the natural action of $(A/lA)^*$ on $Gal(L_{l,z}/L_l)$ to an action of $A/lA$ on $Gal(L_{l,z}/L_l)$. This action is well-defined. Composing with the canonical map from $A$ to $A/lA$, we have an $A$-action on $H_{l,z}$ as well as on $H_{l,\Gamma}$.

The above definition obviously gives that $a \cdot \psi_\lambda = \psi_{\phi_a(\lambda)}$ for $a \in A$. In particular, if $\lambda \in H_l$, then so is $\phi_a(\lambda)$ for all $a \in A$. Thus, $H_l$ is an $A$-submodule of $\Lambda_l^\phi$. To summarize the above discussion, we have the following results as in Proposition 2.3:

**Proposition 3.1.** *Let $l$ be a monic irreducible polynomial in $A$ satisfying the above conditions. We have:*

(1). *The $A$-action defined as above gives an $A$-module structure on $H_{l,z}$ and consequently gives an $A$-module structure on $H_{l,\Gamma}$.*

(2). $H_l$ is an $A$-submodule of $\Lambda_l^\phi$ and $H_{l,z}, H_l$ are isomorphic as $A$-modules. Consequently, $H_{l,z}$ and $H_{l,\Gamma}$ are $A$-modules of exponent $l$.

Let $\Gamma$ be a finitely generated $A$-submodule of the additive group $(L, +)$. Let $\frac{1}{l}\Gamma = \{\alpha \in \bar{k} | \phi_l(\alpha) \in \Gamma\}$ be the $l$-division module of $\Gamma$. Denote by $L_l = L(\Lambda_l^\phi)$ and $L_{l,\Gamma} = L(\Lambda_l^\phi, \frac{1}{l}\Gamma)$. Analogous to the classical Kummer theory, we are interested in the following tower of Kummer extensions $L \subset L_l \subset L_{l,\Gamma}$ with associated Galois groups:

$$\begin{array}{c} L_{l,\Gamma} \\ H_\Gamma(l) \mid \\ L_l \; G_\Gamma(l) \\ G_\Gamma(l)/H_\Gamma(l) \simeq G(l) \mid \\ L \end{array}$$

Since $G(l) \cong (A/lA)^*$ has order prime to the order of $\Lambda_l^\phi$, by a well-known result in [2, Cor. 10.2], we have the following:

**Proposition 3.2.** $H^1(G(l), \Lambda_l^\phi) = 0$.

By (1.6), the Galois group $G(l)$ is isomorphic to $(A/lA)^*$. The main goal is to show that under some mild condition, $H_\Gamma(l)$ is as large as possible.

Let $\Gamma' = \frac{1}{l}\Gamma \cap L$ and define the exponent $e(\Gamma'/\Gamma)$ to be the unique monic polynomial with smallest degree such that $\phi_{e(\Gamma'/\Gamma)}(\Gamma') \subseteq \Gamma$. It is easy to check that $e(\Gamma'/\Gamma)$ is well-defined.

For each $a \in \Gamma$, let $\alpha \in \bar{k}$ be a root of the polynomial $f(X) = \phi_l(X) - a$. Let $\sigma \in H_{l,\Gamma}$. Then $\sigma\alpha = \alpha + \lambda_\sigma$ for some $\lambda_\sigma \in \Lambda_l^\phi$. The map $\sigma \mapsto \lambda_\sigma$ is obviously a homomorphism of $H_{l,\Gamma}$ into $\Lambda_l^\phi$. Write $\lambda_\sigma = \sigma\alpha - \alpha$. It is easy to see that $\lambda_\sigma$ is independent of the choice of the root $\alpha$ of $\phi_l(X) - a$. We denote $\lambda_\sigma$ by $<\sigma, a>$. The map $(\sigma, a) \mapsto <\sigma, a>$ gives us a map $H_{l,\Gamma} \times \Gamma \to \Lambda_l^\phi$.

**Proposition 3.3.** The map $H_{l,\Gamma} \times \Gamma \to \Lambda_l^\phi$ given by $(\sigma, a) \mapsto <\sigma, a>$ is $A$-bilinear, so that the kernel on the left is $\{1\}$ and the kernel on the right is $\Gamma \cap \phi_l(L_l)$.

*Proof.* Similar to the proof of Proposition 2.8. □

Thus, we have an $A$-module homomorphism $\varphi : \Gamma \to Hom_A(H_\Gamma(l), \Lambda_l^\phi)$. More precisely, for each $a \in \Gamma$, we have an $A$-module homomorphism

$$\varphi_a : H_\Gamma(l) \to \Lambda_l^\phi \text{ defined by } \varphi_a(\sigma) = \sigma\alpha - \alpha,$$

where $\phi_l(\alpha) = a$.

By the same way as discussed in Section 2, we can get the following results:

**Theorem 3.4.** Let $e_l(\Gamma) = g.c.d.(e(\Gamma'/\Gamma), l)$ and let $\Gamma_\varphi$ be the kernel of $\varphi$. Then $\phi_{e_l(\Gamma)}(\Gamma_\varphi) \subseteq \phi_l(\Gamma)$.

**Corollary 3.5.** If $e_l(\Gamma) = 1$, i.e. $l \nmid e(\Gamma'/\Gamma)$, then $\phi_l(\Gamma) = \Gamma \cap \phi_l(L_l) = \Gamma \cap \phi_l(L)$. In this case, the pairing $H_\Gamma(l) \times \Gamma/\phi_l(\Gamma) \to \Lambda_l^\phi$ is nondegenerate. Consequently, we have an $A$-module ( resp. $A/lA$-module) isomorphism

$$\varphi : \Gamma/\phi_l(\Gamma) \to Hom_A(H_\Gamma(l), \Lambda_l^\phi) (\text{ resp. } \varphi : \Gamma/\phi_l(\Gamma) \to Hom_{A/lA}(H_\Gamma(l), \Lambda_l^\phi)).$$

**Corollary 3.6.** *If $e_l(\Gamma) = 1$ and $\Gamma$ is free of rank $r$ with basis $\{a_1, ..., a_r\}$, let $\varphi_i = \varphi_{a_i}$, then the map $H_\Gamma(l) \to \Lambda_l^\phi \times ... \times \Lambda_l^\phi$ ( r-copies) given by $\sigma \mapsto (\varphi_1(\sigma), ..., \varphi_r(\sigma))$ is an $A$-module ( resp. $A/lA$-module) isomorphism.*

Note that it is not clear whether $e_l(\Gamma) = 1$ for almost all monic irreducible polynomials $l$ in $A$. In order to obtain the result that $H_\Gamma(l) \cong \Lambda_l^\phi \times ... \times \Lambda_l^\phi$ ( r-copies) for almost all primes $l$ in $A$, we give the proof as follows:

First, recall the following well-known result (see also [12, P.71, Lemma]):

**Lemma 3.7.** *Let $R$ be a product of fields, and let $V$ be a free rank $1$ module over $R$. Suppose that $C$ is an $R$-submodule of $B = V \times ... \times V$ (n times) which is strictly smaller than $B$. Then there are elements $t_1, ..., t_n$ of $R$, not all $0$, such that $\sum t_i v_i = 0$ for all $(v_1, ..., v_n) \in C$.*

By taking $R = A/lA$, $V = \Lambda_l^\phi$ and $C = H_\Gamma(l)$, then it is sufficient to show that there are elements $\varphi_1, ..., \varphi_r \in H_\Gamma(l)$ which are linearly independent over $A/lA$.

Let $H_l = Gal(L^{sep}/L(\Lambda_l^\phi))$. Consider the map $\varphi': L \to Hom(H_l, \Lambda_l^\phi)$ given by $x \mapsto \varphi'_x$, where $\varphi'_x(\sigma) = \sigma(\alpha) - \alpha$ for $\sigma \in H_l$ and some $\alpha$ with $\phi_l(\alpha) = x$. It is easy to see that the map $\varphi'$ is $A$-linear. Consider the map $\delta: L \to H^1(Gal(L^{sep}/L), \Lambda_l^\phi)$, which is obtained by taking cohomology in the short exact sequence $0 \to \Lambda_l^\phi \to L^{sep} \xrightarrow{l} L^{sep} \to 0$. By definition, it is easy to see that $\varphi'$ is the composition of $\delta$ with the restriction homomorphism $Res.: H^1(Gal(L^{sep}/L), \Lambda_l^\phi) \to H^1(H_l, \Lambda_l^\phi) = Hom(H_l, \Lambda_l^\phi)$. By the restriction-inflation sequence together with the vanishing of $H^1(G(l), \Lambda_l^\phi)$ (given in Proposition 3.2), we have that $\varphi'$ induces an $A/lA$-linear injection $L/\phi_l(L) \to Hom(H_l, \Lambda_l^\phi)$. Notice that if we restrict $\varphi'$ to $\Gamma/\phi_l(\Gamma)$, then each $\varphi'_x$ in $\varphi'(\Gamma/\phi_l(\Gamma))$ factors through $H_\Gamma(l)$. So, we may view the map $\varphi'|_{\Gamma/\phi_l(\Gamma)}$ as the natural map $\Gamma/\phi_l(\Gamma) \to Hom(H_\Gamma(l), \Lambda_l^\phi)$ given by $a \mapsto \varphi_a$ as defined above Theorem 3.4. By the same arguments as in [4, Theorem 5], we have that for almost all $l$ in $A$, $a_1, ..., a_r$ are linearly independent modulo $\phi_l(L)$. Hence $\varphi_1, ..., \varphi_r$ are linearly independent over $A/lA$.

**Remark:** Since the orders of $\Lambda_l^\phi$ and $G(l)$ are relatively prime, by [2, Cor. 10.2], we have that $H^2(G(l), \Lambda_l^\phi) = 0$. In this case, the orders of $H_\Gamma(l)$ and $G(l)$ are also relatively prime, so $H^2(G(l), H_\Gamma(l)) = 1$, where $G(l)$ acts on $H_\Gamma(l)$ by conjugation. In particular, the exact sequence $1 \to H_\Gamma(l) \to G_\Gamma(l) \to G(l) \to 1$ is split and $G_\Gamma(l)$ is a semidirect product of $H_\Gamma(l)$ by $G(l)$ ( see [2, Ch.IV]).

**Acknowledgements:**

The authors wish to thank Professor Chih-Nung Hsu for his helpful discussion. Also, we wish to thank the referee for invaluable comments and suggestions.

[†]Department of Mathematics, National Taiwan Normal University, Taipei, Taiwan, Republic of China.
*E-mail address*: `wchi@math.ntnu.edu.tw`

[‡]Department of Mathematics, Fu-Jen University, Taipei, Taiwan, Republic of China.
*E-mail address*: `anlyli@math.fju.edu.tw`